\documentclass{amsart}
\numberwithin{equation}{section}

\newtheorem{thm}{Theorem}[section]
\newtheorem{lem}[thm]{Lemma}

\newtheorem{defin}[thm]{Definition}
\newtheorem{rem}[thm]{Remark}

\begin{document}
\title{Initial-boundary value problem}

\author{Ravshan Ashurov}
\author{Oqila Muhiddinova}
\address{National University of Uzbekistan named after Mirzo Ulugbek and Institute of Mathematics, Uzbekistan Academy of Science}
\curraddr{Institute of Mathematics, Uzbekistan Academy of Science,
Tashkent, 81 Mirzo Ulugbek str. 100170} \email{ashurovr@gmail.com}

\small

\title[Initial-boundary value problem]
{Initial-boundary value problem for a time-fractional subdiffusion
equation with an arbitrary elliptic differential operator}

\begin{abstract}

An initial-boundary value problem for a time-fractional
subdiffusion equation with an arbitrary order elliptic
differential operator is considered. Uniqueness and existence of
the classical solution of the posed problem are proved by the
classical Fourier method. Sufficient conditions for the initial
function and for the right-hand side of the equation are
indicated, under which the corresponding Fourier series converge
absolutely and uniformly. In the case of an initial-boundary value
problem on $N$-dimensional torus, one can easily see that these
conditions are not only sufficient, but also necessary.

\end{abstract}+-

\maketitle

\section{Main results}

The theory of differential equations with fractional derivatives
has gained considerable popularity and importance in the past few
decades, mainly due to its applications in numerous seemingly
distant fields of science and technology (see, for example,
\cite{Mach} - \cite{Hil}). As the most important applications of
this theory one may consider the recent investigations on modeling
of COVID-19 outbreak \cite{YChen}, \cite{Khan}. The data \cite{JH}
presented by Johns Hopkins University on the outbreak from
different countries seem to show fractional order dynamical
processes.

In tern, the well-deserved popularity of the theory attracted the
attention of specialists, causing a large number of investigations
on mathematical aspects of fractional differential equations and
methods to solve them (see, for example, \cite{Mach} and
references therein, \cite{SU} - \cite{ACT}).

In this paper we will investigate by the Fourier method the
solvability (in the classical sense) of initial-boundary value
problems for a time-fractional subdiffusion equation with elliptic
differential operators of any order, defined on an arbitrary
$N$-dimensional bounded domain $\Omega$ with smooth boundary
$\partial\Omega$. The fractional part of our equation will be the
Riemann - Liouville fractional derivative $\partial_t^\rho$ of
order $0<\rho\leq 1$.

Many special cases of this problem have been considered by a
number of authors using different methods. It has been mainly
considered  the case of one spatial variable $x\in \mathbb{R}$ and
subdiffusion equation with "elliptical part" $u_{xx}$ (see, for
example, handbook \cite{Mach}, book of A.A. Kilbas et al.
\cite{Kil} and monograph of A. V. Pskhu \cite{PSK}, and references
in these works). In multidimensional case ($x\in \mathbb{R}^N$)
instead of the differential expression $u_{xx}$ it has been
considered either the Laplace operator (\cite{Kil}, \cite{PS2}),
or the elliptic differential or pseudo-differential operator in
the whole space $\mathbb{R}^N$ with constant coefficients
\cite{SU}. In both cases the authors investigated the Cauchy type
problems applying either the Laplace transform or the Fourier
transform. In his recent paper \cite{PS1} A. V. Pskhu considered
an initial-boundary value problem for subdiffusion equation with
the Laplace operator and domain $\Omega$ - a multidimensional
rectangular region. The author succeeded to construct the Green's
function. Naturally, based on the physical meaning of the problem,
it is of interest to consider similar problems for arbitrary
bounded spatial domains $\Omega$.

In an arbitrary domain $\Omega$ initial-boundary value problems
for subdiffusion equations (the fractional part of the equation is
a multi-term and initial conditions are non-local) with the Caputo
derivatives has been investigated by M. Ruzhansky et al.
\cite{Ruz}. The authors proved the existence and uniqueness of the
generalized solution to the problem.

Let $A(x,D)=\sum\limits_{|\alpha|\leq m} a_\alpha (x) D^\alpha $
be an arbitrary positive formally selfadjoint (symmetric) elliptic
differential operator of order $m=2l$ with sufficiently smooth
coefficients $a_\alpha (x)$ in $\Omega$, where $\alpha=(\alpha_1,
\alpha_2, ..., \alpha_N)$ - multi-index and $D=(D_1, D_2, ...,
D_N)$, $D_j=\frac{\partial}{\partial x_j}$. Recall, an operator
$A(x,D)$ is elliptic in $\Omega$, if for all $x\in \Omega$ and
$\xi\in \mathbb{R}^N$ one has
$$
\sum\limits_{|\alpha|= m} a_\alpha (x) \xi^\alpha>0 \quad \xi \neq
0.
$$

The fractional integration of order $\rho<0$ of a function $f$
defined on $[0, \infty)$ in the Riemann - Liouville sense is
defined by the formula
$$
\partial_t^\rho f(t)=\frac{1}{\Gamma
(-\rho)}\int\limits_0^t\frac{f(\xi)}{(t-\xi)^{\rho+1}} d\xi, \quad
t>0,
$$
provided the right-hand side exists. Here $\Gamma(\rho)$ is
Euler's gamma function.  Using this definition one can define the
Riemann - Liouville fractional derivative of order $\rho$,
$k-1<\rho\leq k$, $k\in \mathbb{N}$, as (see, for example,
\cite{PSK}, p. 14)
$$
\partial_t^\rho f(t)= \frac{d^k}{dt^k}\partial_t^{\rho-k} f(t).
$$
Note if $\rho=k$, then fractional derivative coincides with the
ordinary classical derivation: $\partial_t^k f(t)
=\frac{d^k}{dt^k} f(t)$.

Let $\rho\in (0, 1]$ be a constant number. Consider the
differential equation
\begin{equation}\label{eq}
\partial_t^\rho u(x,t) + A(x,D)u(x,t) = f(x,t), \quad x\in \Omega, \quad t>0;
\end{equation}
with initial
\begin{equation}\label{in}
\lim\limits_{t\rightarrow 0}\partial_t^{\rho-1} u(x,t) =
\varphi(x), \quad x\in \Omega;
\end{equation}
and boundary
\begin{equation}\label{bo}
B_j u(x,t) =\sum\limits_{|\alpha|\leq m_j} b_{\alpha,j} (x)
D^\alpha u(x,t) = 0, \quad 0\leq m_j\leq m-1, \quad j=1,2,..., l;
\quad x\in
\partial \Omega;
\end{equation}
conditions, where $f(x,t), \varphi(x)$ and coefficients $
b_{\alpha,j} (x)$ are given functions.

\begin{defin} A function $u(x,t)$ with the properties
$\partial_t^\rho u(x,t), A(x,D)u(x,t)\in C(\bar{\Omega}\times
(0,\infty))$, $\partial_t^{\rho-1} u(x,t)\in C(\bar{\Omega}\times
[0,\infty))$  and satisfying all the conditions of problem
(\ref{eq}) - (\ref{bo}) in the classical sense is called
\textbf{the regular solution} of initial-boundary value problem
(\ref{eq}) - (\ref{bo}).
\end{defin}

We draw attention to the fact, that in the definition of the
regular solution the requirement of continuity in the closed
domain of all derivatives included in equation (\ref{eq}) is not
caused by the merits. However, on the one hand, the uniqueness of
just such a solution is proved quite simply, and on the other, the
solution found by the Fourier method satisfies the above
conditions.

We will also call the regular solution simply the solution to the
boundary value problem.

Application of the Fourier method to problem (\ref{eq}) -
(\ref{bo}) leads us to consider the following spectral problem
\begin{equation}\label{sur}
A(x,D)v(x) =\lambda v(x) \quad x\in \Omega;
\end{equation}
\begin{equation}\label{sgu}
B_j v(x) = 0, \quad j=1,2,..., l; \quad x\in \partial \Omega.
\end{equation}

S. Agmon \cite{Agm} found the necessary conditions for boundary
$\partial\Omega$ of the domain $\Omega$ and for coefficients of
operators $A$ and $B_j$, which guarantee compactness of the
inverse operator, i.e. the existence of a complete in
$L_2(\Omega)$ system of orthonormal eigenfunctions $\{v_k(x)\}$
and a countable set of positive eigenvalues $\lambda_k$ of
spectral problem (\ref{sur}) - (\ref{sgu}). We will call these
conditions as \emph{condition} $(A)$.

In accordance with the Fourier method, we will look for a solution
to problem  (\ref{eq}) - (\ref{bo}) in the form of a series:
$$
u(x,t) = \sum\limits_1^\infty T_j(t) v_j(x),
$$
where functions $T_j(t)$ are solutions of the Cauchy type problem
\begin{equation}\label{ca}
\partial_t^\rho T_j + \lambda_j T_j =f_j(t),\quad \lim\limits_{t\rightarrow 0}\partial_t^{\rho-1}
T_j(t)=\varphi_j.
\end{equation}
Here we denoted by $f_j(t)$ and $\varphi_j$ the Fourier
coefficients of functions $f(x,t)$ and $\varphi(x)$ with respect
to the system of eigenfunctions $\{v_k(x)\}$ correspondingly,
defined as a scalar product on $L_2(\Omega)$, i.e. for example,
$\varphi_j=(\varphi, v_j)$. The unique solution of problem
(\ref{ca}) has the form (see, for example,  \cite{PSK}, p. 16)
\begin{equation}\label{se2}
T_j(t)=\varphi_j t^{\rho-1} E_{\rho, \rho} (-\lambda_j
t^\rho)+\int\limits_0^t f_j(t-\xi)\xi^{\rho-1} E_{\rho,
\rho}(-\lambda_j\xi^\rho) d\xi,
\end{equation}
where $E_{\rho, \mu}$ is the Mittag-Leffler function of the form
$$
E_{\rho, \mu}(t)= \sum\limits_{k=0}^\infty \frac{t^k}{\Gamma(\rho
k+\mu)}.
$$

A uniqueness result for the solution of problem (\ref{eq}) -
(\ref{bo}) can be formulated as

\begin{thm}\label{uniq} Let condition $(A)$ be satisfied and  $f(x,t) \in C(\bar{\Omega}\times
(0,\infty))$ and $\varphi\in C(\bar{\Omega})$. Then problem
(\ref{eq}) - (\ref{bo}) may have only one regular solution.
\end{thm}

Throughout what follows, we assume that the condition $ (A) $ is
satisfied.

To formulate the existence theorem we need to introduce for any
real number $\tau$ an operator $\hat{A}^\tau$, acting in
$L_2(\Omega)$ in the following way
$$
\hat{A}^\tau g(x)= \sum\limits_{k=1}^\infty \lambda_k^\tau g_k
v_k(x), \quad g_k=(g,v_k).
$$
Obviously, operator $\hat{A}^\tau$ with the domain of definition
$$
D(\hat{A}^\tau)=\{g\in L_2(\Omega):  \sum\limits_{k=1}^\infty
\lambda_k^{2\tau} |g_k|^2 < \infty\}
$$
is selfadjoint. If we denote by $A$ the operator in $L_2(\Omega)$,
acting as $Ag(x)=A(x,D) g (x)$ and with the domain of definition
$D(A)=\{g\in C^m(\bar{\Omega}) : B_j g(x)=0, \quad j=1, ..., l,
\quad x\in
\partial\Omega\}$, then it is not hard to show, that operator $\hat{A}\equiv\hat{A}^1$
is a selfadjoint extension in $L_2(\Omega)$ of operator $A$.

\begin{thm}\label{exis} Let $\tau > \frac{N}{2m}$ and $\varphi\in D(\hat{A}^{\tau})$. Moreover, let $f(x,t)
\in D(\hat{A}^{\tau})$ for all $t\in [0, T]$, and the function
$$
F(t)=\hat{A}^{\tau}f(x,t)
$$
be continuous in the norm of space $L_2(\Omega)$ for all $t\in [0,
T]$. Then there exists a solution of initial-boundary value
problem (\ref{eq}) - (\ref{bo}) and it has the form of series
\begin{equation}\label{se1}
u(x,t) = \sum\limits_1^\infty \big[\varphi_j t^{\rho-1} E_{\rho,
\rho} (-\lambda_j t^\rho)+\int\limits_0^t f_j(t-\xi)\xi^{\rho-1}
E_{\rho, \rho}(-\lambda_j\xi^\rho) d\xi\big] v_j(x),
\end{equation}
which absolutely and uniformly converges on $x\in \bar{\Omega}$
and for each $t\in (0, T]$.

\end{thm}

\begin{rem}\label{cl} Note, when $\rho =1$ this theorem states existence
of the unique solution of the problem, when the equation and
initial condition have the form
$$
u_t(x,t)+A(x,D)u(x,t)=f(x,t),\quad x\in\Omega,\quad t>0;\quad
u(x,0)=\varphi(x),\quad x\in\Omega,
$$
and the boundary condition is the same as (\ref{bo}).
\end{rem}

It should be noted that Sh.A. Alimov in his paper  \cite{Ali}
presents sufficient conditions for a given function to belong to
the domain of definition of operator $\hat{A}^\tau$ in terms of
various classes of differentiable functions. At the end of this
paper we consider initial-boundary value problem
(\ref{eq})-(\ref{bo}) in torus $\mathbb{T}^N$ when operator
$A(x,D)$ has constant coefficients. We will see that in this case
the domain of definition $D(\hat{A}^\tau)$ coincides with the
corresponding Sobolev spaces.

\section{Uniqueness}

In this section we prove Theorem \ref{uniq}. Suppose that all the
conditions of this theorem are satisfied and let initial-boundary
value problem (\ref{eq}) - (\ref{bo}) have two regular solutions
$u_1(x,t)$ and $u_2(x,t)$. Our aim is to prove that
$u(x,t)=u_1(x,t)-u_2(x,t)\equiv 0$. Since the problem is linear,
then we have the following homogenous problem for $u(x,t)$:
\begin{equation}\label{ur1}
\partial_t^\rho u(x,t) + A(x,D)u(x,t) = 0, \quad x\in \Omega, \quad t>0;
\end{equation}
\begin{equation}\label{nu1}
\lim\limits_{t\rightarrow 0} \partial_t^{\rho-1}u(x,t) = 0, \quad
x\in \Omega;
\end{equation}
\begin{equation}\label{gu1}
B_j u(x,t) =\sum\limits_{|\alpha|\leq m_j} b_{\alpha,j} (x)
D^\alpha u(x,t) = 0, \quad 0\leq m_j\leq m-1, \quad j=1,2,..., l;
\quad x\in
\partial \Omega;
\end{equation}

Let $u(x,t)$ be a regular solution of problem
(\ref{ur1})-(\ref{gu1}) and $v_k$ be an arbitrary eigenfunction of
the problem (\ref{sur})-(\ref{sgu}) with the corresponding
eigenvalue $\lambda_k$. Consider the function
\begin{equation}\label{w}
w_k(t)=\int\limits_\Omega u(x,t)v_k(x)dx.
\end{equation}
By definition of the regular solution we may write
$$
\partial_t^\rho w_k(t)=\int\limits_\Omega \partial_t^\rho u(x,t)v_k(x)dx=
-\int\limits_\Omega A(x, D)u(x,t)v_k(x)dx, \quad t>0,
$$
or, integrating by parts,
$$
\partial_t^\rho w_k(t)=-\int\limits_\Omega u(x,t)A(x, D)v_k(x)dx=-\lambda_k
\int\limits_\Omega u(x,t)v_k(x)dx = -\lambda_k w_k(t), \quad t>0.
$$
Since condition (\ref{in}) can be rewritten as (see, for example,
\cite{PSK} p. 104)
\begin{equation}\label{in1}
\lim\limits_{t\rightarrow
0}t^{1-\rho}u(x,t)=\frac{\varphi(x)}{\Gamma(\rho)},
\end{equation}
then, using in (\ref{w}) the homogenous initial condition
(\ref{nu1}) in this form, we have the following Cauchy problem for
$w_k(t)$:
$$
\partial_t^\rho w_k(t) +\lambda_k w_k(t)=0,\quad t>0; \quad
\lim\limits_{t\rightarrow 0}\partial_t^{\rho-1} w_k(t)=0.
$$
This problem has the unique solution; therefore, the function
defined by (\ref{w}), is identically zero: $w_k(t)\equiv 0$ (see
(\ref{se2})). From completeness in $L_2(\Omega)$  of the system of
eigenfunctions $\{v_k(x)\}$, we have $u(x,t) = 0$ for all $x\in
\bar{\Omega}$ and $t>0$. Hence Theorem \ref{uniq} is proved.

\section{Existence}

Here we have borrowed some original ideas from the method
developed in the work of M.A. Krasnoselski et al. \cite{Kra}. The
following lemma plays a key role in this method (\cite{Kra},
p.453).

\begin{lem}\label{CL} Let $\sigma > 1+\frac{N}{2m}$. Then for any $|\alpha|\leq m$
operator $D^\alpha \hat{A}^{-\sigma}$ (completely) continuously
maps from $L_2(\Omega)$ into $C(\bar{\Omega})$ and moreover the
following estimate holds true
\begin{equation}\label{CL}
||D^\alpha \hat{A}^{-\sigma} g||_{C(\Omega)} \leq C
||g||_{L_2(\Omega)}.
\end{equation}

\end{lem}

Using this lemma, we prove that one can validly apply  the
operators $D^\alpha$ with $|\alpha|\leq m$ and $\partial^\rho_t$
to the series (\ref{se1}) term-by-term.

For the Mittag-Leffler function with a negative argument we have
an estimate $|E_{\rho,\mu}(-t)|\leq \frac{C}{1+ t}$ (see, for
example, \cite{PSK},  p.13). Therefore, since all eigenvalues
$\lambda_j$ are positive,
\begin{equation}\label{m1}
|t^{\rho-1} E_{\rho,\rho}(-\lambda_j t^\rho)|\leq
\frac{Ct^{\rho-1}}{1+\lambda_jt^\rho}\leq \frac{C}{\lambda_j
t}(t^\rho\lambda_j)^{\varepsilon/\rho}, \quad t>0,
\end{equation}
where $0<\varepsilon<\rho$. Indeed, let $t^\rho\lambda_j<1$, then
$$
\frac{1}{\lambda_j
t}(t^\rho\lambda_j)^{\varepsilon/\rho}>\frac{1}{\lambda_j
t}t^\rho\lambda_j>t^{\rho-1},
$$
and if $t^\rho\lambda_j>1$, then
$$
\frac{1}{\lambda_j
t}(t^\rho\lambda_j)^{\varepsilon/\rho}>\frac{1}{\lambda_j t}.
$$

Note the series (\ref{se1}) is in fact the sum of two series.
Consider the first series:
\begin{equation}\label{S}
S^1_k(x,t)=\sum\limits_{j=1}^k v_j(x)\varphi_jt^{\rho-1}
E_{\rho,\rho}(-\lambda_j t^\rho),
\end{equation}
and suppose that function $\varphi$ satisfies the condition of
Theorem \ref{exis}, i.e. for some $\tau> \frac{N}{2m}$
$$
\sum\limits_1^\infty \lambda_j^{2\tau} |\varphi_j|^2 \leq
C_\varphi<\infty.
$$
We choose a small $\varepsilon>0$ in such a way, that
$\tau+1-\varepsilon/\rho> 1+\frac{N}{2m}$.  Since
$\hat{A}^{-\tau-1+\varepsilon/\rho} v_j(x) =
\lambda_j^{-\tau-1+\varepsilon/\rho} v_j(x)$, we may rewrite the
sum (\ref{S}) as
$$
S^1_k(x,t)=\hat{A}^{-\tau-1+\varepsilon/\rho}\sum\limits_{j=1}^k
v_j(x)\lambda_j^{\tau+1-\varepsilon/\rho}\varphi_jt^{\rho-1}
E_{\rho,\rho}(-\lambda_j t^\rho).
$$
Therefore by virtue of Lemma \ref{CL} one has
$$
||D^\alpha S^1_k||_{C(\Omega)}=||D^\alpha
\hat{A}^{-\tau-1+\varepsilon/\rho}\sum\limits_{j=1}^k
v_j(x)\lambda_j^{\tau+1-\varepsilon/\rho}\varphi_jt^{\rho-1}
E_{\rho,\rho}(-\lambda_j t^\rho)||_{C(\Omega)}\leq
$$
\begin{equation}\label{S1}
\leq C ||\sum\limits_{j=1}^k
v_j(x)\lambda_j^{\tau+1-\varepsilon/\rho}\varphi_jt^{\rho-1}
E_{\rho,\rho}(-\lambda_j t^\rho)||_{L_2(\Omega)}.
\end{equation}
Using the orthonormality of the system $\{v_j\}$, we will have
\begin{equation}\label{S2}
||D^\alpha S^1_k||^2_{C(\Omega)}\leq C \sum\limits_{j=1}^k
|\lambda_j^{\tau+1-\varepsilon/\rho}\varphi_jt^{\rho-1}
E_{\rho,\rho}(-\lambda_j t^\rho)|^2.
\end{equation}
Application of inequality (\ref{m1}) gives
$$
\sum\limits_{j=1}^k
|\lambda_j^{\tau+1-\varepsilon/\rho}\varphi_jt^{\rho-1}
E_{\rho,\rho}(-\lambda_j t^\rho)|^2\leq C
t^{2(\varepsilon-1)}\sum\limits_{j=1}^k
\lambda_j^{2\tau}|\varphi_j|^2\leq C
t^{2(\varepsilon-1)}C_\varphi.
$$
Therefore we can rewrite the estimate (\ref{S2}) as
$$
||D^\alpha S^1_k||^2_{C(\Omega)}\leq C t^{2(\varepsilon-1)}
C_\varphi.
$$

This implies uniformly on $x\in\bar{\Omega}$ convergence of the
differentiated sum (\ref{S}) with respect to the variables $x_j$
for each $t\in (0,T]$. On the other hand, the sum (\ref{S1})
converges for any permutation of its members as well, since these
terms are mutually orthogonal. This implies the absolute
convergence of the differentiated sum (\ref{S}) on the same
interval $t\in (0,T]$.

Now we consider the second part of the series (\ref{se1}):
\begin{equation}\label{S3}
S^2_k(x,t)=\sum\limits_{j=1}^k v_j(x)\int\limits_0^t
f_j(t-\xi)\xi^{\rho-1} E_{\rho,\rho}(-\lambda_j \xi^\rho)d\xi
\end{equation}
and suppose that function  $f(x,t)$  satisfies all the conditions
of Theorem \ref{exis}, i.e. the following series converges
uniformly on $t\in [0, T]$ for some $\tau> \frac{N}{2m}$:
$$
\sum\limits_1^\infty \lambda_j^{2\tau} |f_j(t)|^2 \leq C_f<\infty.
$$
We have
$$
S^2_k(x,t)=\hat{A}^{-\tau-1+\varepsilon/\rho}\sum\limits_{j=1}^k
v_j(x)\int\limits_0^t
\lambda_j^{\tau+1-\varepsilon/\rho}f_j(t-\xi)\xi^{\rho-1}
E_{\rho,\rho}(-\lambda_j \xi^\rho)d\xi.
$$
Then by virtue of Lemma \ref{CL} one has
$$
||D^\alpha S^2_k||_{C(\Omega)}=||D^\alpha
\hat{A}^{-\tau-1+\varepsilon/\rho}\sum\limits_{j=1}^k
v_j(x)\int\limits_0^t
\lambda_j^{\tau+1-\varepsilon/\rho}f_j(t-\xi)\xi^{\rho-1}
E_{\rho,\rho}(-\lambda_j \xi^\rho)d\xi||_{C(\Omega)}\leq
$$
\begin{equation}\label{S4}
\leq C ||\sum\limits_{j=1}^k v_j(x)\int\limits_0^t
\lambda_j^{\tau+1-\varepsilon/\rho}f_j(t-\xi)\xi^{\rho-1}
E_{\rho,\rho}(-\lambda_j \xi^\rho)d\xi||_{L_2(\Omega)}.
\end{equation}
Using the orthonormality of the system $\{v_j\}$, we will have
$$
||D^\alpha S^2_k||^2_{C(\Omega)}\leq C \sum\limits_{j=1}^k
|\int\limits_0^t
\lambda_j^{\tau+1-\varepsilon/\rho}f_j(t-\xi)\xi^{\rho-1}
E_{\rho,\rho}(-\lambda_j \xi^\rho)d\xi|^2.
$$
Now we use estimate (\ref{m1}) and apply the generalized Minkowski
inequality. Then
$$
||D^\alpha S^2_k||^2_{C(\Omega)}\leq C(\int\limits_0^t
\xi^{\varepsilon-1} (\sum\limits_{j=1}^k
\lambda_j^{2\tau}|f_j(t-\xi)|^2)^{1/2} d\xi)^2\leq C\cdot C_f\cdot
(\frac{T^\varepsilon}{\varepsilon})^2.
$$
Hence, using the same argument as above, we see that the
differentiated sum (\ref{S3}) with respect to the variables $x_j$
converges absolutely and uniformly on $(x,t)\in \bar{\Omega}\times
[0,T]$.

It is not hard to see that
$$
\partial_t^\rho\sum\limits_1^k T_j(t)v_j(x)=  -\sum\limits_1^k \lambda_j
T_j(t)v_j(x)+\sum\limits_1^k f_j(t)v_j(x)=
$$
$$
-A(x,D)\hat{A}^{-\tau-1+\varepsilon/\rho}\sum\limits_1^k
\lambda^{\tau+1-\varepsilon/\rho}_j
T_j(t)v_j(x)+\hat{A}^{-\tau-1+\varepsilon/\rho}\sum\limits_1^k
\lambda_j^{\tau+1-\varepsilon/\rho}f_j(t)v_j(x).
$$

Absolutely and uniformly convergence of the latter  series can be
proved as above.

Obviously, function (\ref{se1}) satisfies the boundary conditions
(\ref{bo}). Considering the initial condition as in (\ref{in1}),
it is not hard to verify, that this condition is also satisfied.

Thus Theorem \ref{exis} is completely proved.

Observe, taking $\alpha =0$ in the above estimations, one may
obtain an estimate for $||u(x,t)||_{C(\Omega)}$, which gives the
stability of the solution to problem (\ref{eq})-(\ref{bo}).

\section{The initial-boundary value problem on $\mathbb{T}^N$}

In the proof of Theorem \ref{exis} we only use the fact that
elliptic operator $A$ has a complete in $L_2(\Omega)$ set of
orthonormal eigenfunctions and Lemma \ref{CL}. This lemma reduces
the study of uniform convergence to the study of convergence in
$L_2$, where the Parseval equality gives a solution to the problem
immediately.  Therefore, similar to Theorem \ref{exis} statement
holds true for any operator with these properties. As an example
we may consider differential operator with involution and the
Bessel operator (see the paper of M. Ruzhansky et al. \cite{Ruz}),
or the initial-boundary value problem in $N$-dimensional torus:
$\mathbb{T}^N = (\pi, \pi]^N$.

Let us consider the last case. So let
$A(D)=\sum\limits_{|\alpha|=m} a_\alpha D^\alpha$ be a homogeneous
symmetric positive elliptic differential operator with constant
coefficients. Let the differential equation and initial condition
have the form ($0< \rho\leq 1$)
\begin{equation}\label{eqt}
\partial_t^\rho u(x,t) + A(D)u(x,t) = f(x,t), \quad x\in \mathbb{T}^N, \quad t>0;
\end{equation}
\begin{equation}\label{int}
\lim\limits_{t\rightarrow 0}\partial_t^{\rho-1} u(x,t) =
\varphi(x), \quad x\in \mathbb{T}^N;
\end{equation}
Instead of boundary conditions (\ref{bo}) we consider the
$2\pi$-periodic in each argument $x_j$ functions and suppose that
$\varphi(x)$ and $f(x,t)$ are also $2\pi$-periodic functions in
$x_j$.

Let $A$ stands for the operator $A(D)$, defined on $2\pi$-periodic
functions from $C^m(\mathbb{R}^N)$. The closure $\hat{A}$ of this
operator in $L_2(\mathbb{T}^N)$ is selfadjoint. It is not hard to
see that operator $\hat{A}$ has a complete (in
$L_2(\mathbb{T}^N)$) set of eigenfunctions
$\{(2\pi)^{-N/2}e^{inx}\}$, and corresponding eigenvalues $A(n)$.
Therefore, by virtue of J. von Niemann theorem, for any $\tau> 0$
operator $\hat{A}^\tau$ acts as $\hat{A}^\tau
g(x)=\sum\limits_{A(n)<\lambda} A^\tau(n) g_n e^{inx}$, where
$g_n$ is Fourier coefficients of $g\in L_2(\mathbb{T}^N)$. The
domain of definition of this operator is defined from the
condition $\hat{A}^\tau g(x)\in L_2(\mathbb{T}^N)$ and has the
form
$$
D(\hat{A}^\tau)=\{g\in L_2(\mathbb{T}^N): \sum\limits_{n\in
\mathbb{Z}^N} A^{2\tau}(n) |g_n|^2 < \infty\}.
$$

In order to define the domain of definition of operator
$\hat{A}^\tau$ in terms of the Sobolev spaces, we remind the
definition of these spaces (see, for example, \cite{AAP}): we say
that function $g\in L_2(\mathbb{T}^N)$ belongs to the Sobolev
space $L_2^a(\mathbb{T}^N)$ with the real number $a> 0$, if the
norm
$$
||g||^2_{L_2^a(\mathbb{T}^N)}=\big|\big|\sum\limits_{n\in
Z^N}(1+|n|^2)^{\frac{a}{2}}g_n
e^{inx}\big|\big|^2_{L_2(\mathbb{T}^N)}=\sum\limits_{n\in
Z^N}(1+|n|^2)^a|g_n|^2
$$
is bounded. When $a$ is not integer, then this space is also
called the Liouville space.

It is not hard to show the existence of constants $c_1$ and $c_2$
such that
$$
c_1(1+|n|^2)^{\tau m} \leq 1+A^{2\tau}(n)\leq c_2 (1+|n|^2)^{\tau
m}.
$$
Therefore, $ D(\hat{A}^\tau)= L_2^{\tau m}(\mathbb{T}^N)$.

Repeating verbatim the proof of Theorem \ref{exis}, the following
statement is proved.

\begin{thm}\label{TN} Let $\tau > \frac{N}{2}$ and $\varphi\in L^{\tau}_2(\mathbb{T}^N)$. Moreover, let $f(x,t)\in L^\tau_2(\mathbb{T}^N)$
for any $t\in [0,T]$ and the function
$$
F(t)=\hat{A}^{\tau/m}f(x,t)
$$
be continuous in the norm of space $L_2(\mathbb{T}^N)$ for all
$t\in [0, T]$. Then there exists a solution of initial-boundary
value problem (\ref{eqt}) - (\ref{int}) and it has the form
$$
u(x,t)=\sum\limits_{n\in\mathbb{Z}^N} \big[\varphi_n t^{\rho-1}
E_{\rho, \rho} (-A(n) t^\rho)+\int\limits_0^t
f_n(t-\xi)\xi^{\rho-1} E_{\rho, \rho}(-A(n)\xi^\rho)
d\xi\big]e^{inx},
$$
which absolutely and uniformly converges on $x\in \mathbb{T}^N$
and for each $t\in (0, T]$, where $\varphi_n$ and $f_n(t)$ are
corresponding Fourier coefficients.

\end{thm}
\begin{rem}
Note, when $\tau>\frac{N}{2}$,  according to the Sobolev embedding
theorem, all functions in $L^{\tau}_2(\mathbb{T}^N)$ are
$2\pi$-periodic continuous functions. The fulfillment of the
inverse inequality $\tau \leq \frac{N}{2}$, admits the existence
of unbounded functions in $L^{\tau}_2(\mathbb{T}^N)$ (see, for
example, \cite{AAP}). Therefore, condition  $\tau>\frac{N}{2}$ of
this theorem is not only sufficient for the statement to be hold,
but it is also necessary.
\end{rem}

A $2\pi$-periodical continuous function $\varphi(x)$, defined on
$\mathbb{T}^N$, is said to be a H\"{o}lder continuous with
exponent $\beta\in (0,1]$, if for some constant $C$,
$$
|\varphi(x)-\varphi(y)|\leq C|x-y|^\beta,\quad \text{for} \quad
\text{all}\quad x,y\in \mathbb{T}^N.
$$
We define $C^a(\mathbb{T}^N)$ to be the H\"{o}lder space of
$2\pi$-periodical functions in $C^{[a]}(\mathbb{T}^N)$, $[a]$ is
the integer part of $a$, all derivatives $D^\alpha \varphi$,
$|\alpha|=[a]$, of which are H\"{o}lder continuous with exponent
$a-[a]$. It follows from the embedding theorem $C^a(\mathbb{T}^N)
\rightarrow L_2^{a-\varepsilon}(\mathbb{T}^N)$, $\varepsilon>0$,
that the conditions of the Theorem \ref{TN} can be formulated in
terms of the H\"{o}lder spaces, replacing the classes
$L^{\tau}_2(\mathbb{T}^N)$ by the classes $C^\tau(\mathbb{T}^N)$.
Again, condition $\tau>\frac{N}{2}$ is precise: if
$\tau=\frac{N}{2}$, then there exists a function in
$C^\tau(\mathbb{T}^N)$, such that the Fourier series of which
diverges at some point (see, for example, \cite{AAP}).

\section{Acknowledgement} The authors convey thanks to Sh. A.
Alimov and S. R. Umarov for discussions of these results. We  wish
also to thank B. Turmetov for making us  aware  of  the  paper  of
M. Ruzhansky et al. \cite{Ruz}.

\

\

\bibliographystyle{amsplain}

\end{document}